\documentclass[11pt,a4paper]{article}
\usepackage{amsmath,amssymb,amsthm,graphicx,color,bbm,array}
\usepackage[left=1in,right=1in,top=1in, bottom=1in]{geometry}
\usepackage{cite}
\usepackage[british]{babel}
\usepackage{hyperref} 
\usepackage{enumitem}
\usepackage{tikz}
\usepackage{caption}
\usepackage{mathtools}
\usepackage{comment}
\hypersetup{
    colorlinks=false,
    pdfborder={0 0 0},
}

\newtheorem{theorem}{Theorem}[section]
\newtheorem{corollary}[theorem]{Corollary}
\newtheorem{proposition}[theorem]{Proposition}
\newtheorem{lemma}[theorem]{Lemma}

\numberwithin{equation}{section}

\theoremstyle{definition}

\theoremstyle{remark}
\newtheorem{remark}[theorem]{Remark}
\newtheorem{remarks}[theorem]{Remarks}
\newtheorem*{remark*}{Remark}

\newcommand{\bb}[1]{\mathbb{#1}}
\newcommand{\ind}[1]{\mathbbm{1}_{\{#1\}}}

\DeclareMathOperator{\Exp}{\mathbb{E}}
\let\Pr\relax
\DeclareMathOperator{\Pr}{\mathbb{P}}

\newcommand{\R}{\ensuremath{\mathbb{R}}}

\newcommand{\ZP}{\ensuremath{\mathbb{Z}_+}}
\newcommand{\RP}{\ensuremath{\mathbb{R}_+}}
\newcommand{\N}{\ensuremath{\mathbb{N}}}

\newcommand{\ud}{\,\mathrm{d}}
\newcommand{\re}{{\mathrm{e}}}

\newcommand{\eps}{\varepsilon}

\newcommand{\tY}{\widetilde{Y}}
\newcommand{\tA}{\widetilde{A}}

\newcommand{\bV}{\obar{V}}
\newcommand{\cE}{\mathcal{E}}

\newcommand{\as}{\ \text{a.s.}}

\newcommand{\toas}{\overset{\textup{a.s.}}{\longrightarrow}}
\newcommand{\tod}{\overset{\textup{d}}{\longrightarrow}}
\newcommand{\eqd}{\overset{\textup{d}}{=}}

\newcommand{\cF}{\mathcal{F}}

\newcommand{\bessq}{\mathrm{BESQ}}
\newcommand{\besq}[2]{{\bessq}^{#1}(#2)}

\newcommand{\sscbi}[4]{{\mathrm{SSCBI}}({#1},{#2},{#3};{#4})}

\newcommand{\obar}[1]{\mkern 1.5mu\overline{\mkern-1.5mu#1\mkern-1.5mu}\mkern 1.5mu}

\makeatletter
\def\namedlabel#1#2{\begingroup  
    (#2)%
    \def\@currentlabel{#2}%
    \phantomsection\label{#1}\endgroup
}
\makeatother

\newlist{myenumi}{enumerate}{10}
\setlist[myenumi]{leftmargin=0pt, labelindent=\parindent, listparindent=\parindent, labelwidth=0pt, itemindent=!, itemsep=1pt, parsep=4pt}

\newlist{thmenumi}{enumerate}{10}
\setlist[thmenumi]{leftmargin=0pt, labelindent=\parindent, listparindent=\parindent, labelwidth=0pt, itemindent=!}

\title{Superdiffusive limits for stochastic kinetics  \\ driven by self-similar drifts}

\author{Miha Bre\v sar\footnote{\raggedright School of Data Science, The Chinese University of Hong Kong, Shenzhen, 
2001, Longxiang Boulevard, Longgang District, Shenzhen, China. Email:~\href{mailto:mihabresar@cuhk.edu.cn}{\texttt{mihabresar@cuhk.edu.cn}}}
\and Conrado da Costa\footnote{\raggedright Institute of Mathematics and Statistics (IME), University of São Paulo (USP), R. do Matão, 1010 - Butantã, São Paulo, SP,  Brazil. Email:~\href{mailto:conrado.dacosta@ime.usp.br}{\texttt{conrado.dacosta@ime.usp.br}}}
\and Aleksandar Mijatovi\'c\footnote{\raggedright Department of Statistics, University of Warwick, 
Coventry, CV4 7AL, UK. Email:~\href{mailto:a.mijatovic@warwick.ac.uk}{\texttt{a.mijatovic@warwick.ac.uk}}}~\footnote{The Alan Turing Institute.} \and Andrew Wade\footnote{\raggedright Department of Mathematical Sciences, Durham University, Upper Mountjoy, Durham, DH1 3LE, UK. Email:~\href{mailto:andrew.wade@durham.ac.uk}{\texttt{andrew.wade@durham.ac.uk}}.} }
\date{\today}

\begin{document}

\maketitle
\renewcommand{\thefootnote}{\textdagger}

\begin{abstract}
We prove anomalous-diffusion scaling for a one-dimensional stochastic kinetic 
dynamics, in which the stochastic drift is driven by an exogenous self-similar noise, and also includes endogenous volatility which is permitted to have arbitrary dependence with the exogenous noise. We identify the superdiffusive scaling exponent for the model, and prove strong and weak convergence results on the corresponding scale. Our framework admits self-similar noise that is either a Bessel process, or, more generally, a self-similar continuous-state branching process with immigration, as well as more general processes satisfying certain asymptotic conditions. 
\end{abstract}

\medskip


\medskip


\maketitle

\section{Introduction and main results}
\label{sec:results}

The subject of this work is the long-term behaviour of an It\^{o} process $X=(X_t)_{t\in\RP}$ 
in $\RP:=[0,\infty)$
with a representation 
\begin{equation}
\label{eq:X_dynamics}
 X_t = X_0+\int_0^t \frac{f(s,Y_s)}{X_s}\ud s + B_t, \qquad \text{ $X_0>0$,}
\end{equation}
where the  process $Y=(Y_t)_{t\in\RP}$ is adapted to the same filtration as the Brownian motion~(BM) $B=(B_t)_{t\in\RP}$.
Under suitable assumptions on the function 
$f: \RP\times\RP \to \RP$
and the exogenous noise process $Y$, which drives the \emph{stochastic drift} of~$X$,
the process $X$ will exhibit \emph{anomalous diffusion}.
 Physical motivation includes dynamics of particles interacting with an external field or
 medium, or with an internal relaxation mechanism; see e.g.~\cite{bg,FT20,FT21} and~\S\ref{sec:literature} below for further discussion of motivation and related literature. 
 
 We assume the  function~$f(t,y)$ 
 that contributes to the drift via~\eqref{eq:X_dynamics}
 has certain polynomial asymptotic growth 
 behaviour for large $t$ and $y$: informally,
that $f(t,y) \sim \rho t^\gamma y^\alpha$ as both $t, y\to\infty$,
for a constant $\rho>0$. In fact, the hypothesis is a little
stronger than this (giving some comparable control for fixed~$t$ and large~$y$, for example).
 
\renewcommand{\thefootnote}{\arabic{footnote}}
\setcounter{footnote}{0} 
\begin{description}
\item
[\namedlabel{ass:parameters}{A$_f$}]
Suppose that $f: \RP\times\RP \to \RP$ is continuous, and 
that for some constant
$\rho \in (0,\infty)$
and exponents 
$\alpha\in \RP$ and $\gamma\in\R$
satisfies the following.
For every $\eps>0$, there is some $r_\eps \in \RP$ such that 
\begin{equation}
\label{eq:A}
\sup_{(t,y)\in  \RP^2: \, t+ y \geq r_\eps} \left| f (t, y) (1+t)^{-\gamma} (1+y)^{-\alpha} - \rho \right| \leq \eps.
\end{equation}
\end{description}
 
The aim of this paper is to give some natural and robust growth and stability hypotheses on the exogenous noise process~$Y$
so that process~$X$ satisfying~\eqref{eq:X_dynamics} exhibits superdiffusive asymptotic behaviour, as quantified via (i) a distributional scaling limit, and (ii) an a.s.-quantification of the transient growth exponent. Our general result is presented in Theorem~\ref{thm:general} below. 
We first state our principal concrete result, in Theorem~\ref{thm:limit-new}, for SSCBI exogenous noise (see~\eqref{eq:sscbi} below for a definition); the important special case of squared-Bessel noise is then stated as Corollary~\ref{cor:bes-limit}.
 First we need some notation.

A non-negative \emph{self-similar continuous-state branching process with immigration} (SSCBI) is a $[0,\infty]$-valued stochastic process $Y = (Y_t)_{t \in \RP}$ whose law is determined by the Laplace functional
\begin{equation}
    \label{eq:sscbi}
 \Exp_y [ \re^{-\lambda Y_t} ] = (1 + q c_0 \lambda^q t )^{-\delta} \exp \left( - \frac{\lambda y}{(1+q c_0 \lambda^q t)^{1/q}} \right) ,\end{equation}
where $\lambda >0$, $t \geq 0$, and $y \geq 0$, and
parameters $q \in (0,1]$, $c_0>0$, $\delta >0$. 
We write $Y \sim \sscbi{q}{c_0}{\delta}{y}$ if $Y$ has law determined by~\eqref{eq:sscbi} and has initial value $Y_0 = y \in \RP$.
Write $\kappa := 1/q \in [1,\infty)$. The SSCBI  are a subset of  continuous-state branching process with immigration   that arise as weak limits of discrete branching processes with immigration for correctly-scaled offspring and immigration distributions as shown in~\cite{kw}; for a recent introduction to CBI processes see~\cite{li}, and some other modern presentations and results can be found in~\cite{fub,krm,cpu,mpu}. 
The case $q = \kappa =  1$ is special in that $Y$ has continuous sample paths; indeed, when
$c_0 = 2$ it is a squared Bessel process of dimension $2\delta$,
and a multiple of it otherwise.

Note from~\eqref{eq:sscbi} that, for $c>0$, $\Exp_y [ \re^{-\lambda c^{-\kappa} Y_{ct} } ] = \Exp_{c^{-\kappa}y} [ \re^{-\lambda Y_t} ]$,
which demonstrates that $1/\kappa = q$ is the \emph{self-similarity index} of the SSCBI. In particular, it follows that $t^{-\kappa} Y_t$, with $Y$ started from any $y \geq 0$, converges in distribution as $t \to \infty$ to $Y_1$ started from $0$.

The following is our general result as applied to the case of SSCBI exogenous noise. Because of the apparent singularity in the drift~\eqref{eq:X_dynamics}, in the same way as one does for the Bessel process, 
it is more convenient to formulate the dynamics via ``$X$ squared'', which is essentially the process~$S$ in~\eqref{eq:sq-X-dynamics}: see Remark~\ref{rems:limit}\ref{rems:limit-ii} after the statement of the theorem.
We denote convergence in distribution by `$\tod$'.

\begin{theorem}
\label{thm:limit-new} 
Let $f:\RP\times\RP\to\RP$ satisfy~\eqref{ass:parameters},
$\gamma+\alpha \kappa >0$,
and assume $q \in (0,1]$ and $\delta \in (0,\infty)$.
Suppose that either (i) $q=1$; or (ii) $0 \leq \alpha < q$.
Suppose 
the adapted process $(S,Y,B)$ consists of 
an $\R$-valued Brownian motion $B$, 
 $Y \sim \sscbi{q}{c_0}{\delta}{y}$  
with self-similarity index $1/\kappa = q$
 and an $\RP$-valued process
$S=(S_t)_{t\in\RP}$, satisfying 
\begin{equation}
    \label{eq:sq-X-dynamics}
     S_t  = S_0 + \int_0^t \big(2f(s,  Y_s ) + 1\big) \ud s + 2 \int_0^t \sqrt{ {S_s } } \ud B_s, \text{ for all } t \in \RP,
   \end{equation}
   started at a deterministic level $S_0\in\RP$. 
Then 
the $\RP$-valued process $X=(X_t)_{t\in\RP}$,
given by $X_t := \sqrt{ S_t }$, $t \in \RP$, 
has the following asymptotic properties.
\begin{thmenumi}[label=(\alph*)]
\item\label{thm:limit-a2} There is superdiffusive transience, namely, 
\begin{equation}
    \label{eq:escape-rate2}
\lim_{t \to \infty} \frac{\log X_t}{\log t} = \frac{1+\gamma+\alpha\kappa}{2}, \as
\end{equation}
\item\label{thm:limit-b2}
There is a distributional limit, namely, 
\begin{equation}
    \label{eq:limit-statement2}
t^{-(1+\gamma+\alpha\kappa)/2}X_t \tod   \left(2\rho \int_0^1 s^\gamma \tY_s^\alpha \ud s\right)^{1/2},  \text{ as } t \to \infty,
\end{equation}
where the limit is a.s.~finite and $\tY$ in the limit has the same law as $Y$ but started at $0$.
\end{thmenumi}
\end{theorem}

Denote by $\besq{\delta}{y}$  the law of
a squared-Bessel process $Y$ with ``dimension'' parameter $\delta\in(0,\infty)$ started at an arbitrary $Y_0 = y\in\RP$;
this law can be defined as that of the solution of the stochastic differential equation (SDE) in~\eqref{eq:sq_bessel} below.
Recall that
$\delta > 0$ ensures  $Y$   is not absorbed at $0$;
for $\delta \in (0,2)$, the process~$Y$ is point recurrent over~$\RP$, while for $\delta \geq 2$ it is point transient (see~\cite[Ch.~XI]{RevYor04} for further details). Since the SSCBI with $q = \kappa = 1$ and $c_0 =2$ reduces to the squared Bessel process, the following is a direct consequence of Theorem~\ref{thm:limit-new}.

\begin{corollary}
\label{cor:bes-limit} 
Let $f:\RP\times\RP\to\RP$ satisfy~\eqref{ass:parameters}
with $\alpha+\gamma >0$, 
and assume $\delta \in (0,\infty)$.
Suppose 
the adapted process $(S,Y,B)$ consists of 
an $\R$-valued Brownian motion $B$, 
a squared-Bessel process $Y$ with law $\besq{\delta}{y}$ (defined via SDE~\eqref{eq:sq_bessel}) and an $\RP$-valued process
$S=(S_t)_{t\in\RP}$, satisfying 
SDE~\eqref{eq:sq-X-dynamics}
   started at a deterministic level $S_0\in\RP$. 
Then 
the $\RP$-valued process $X=(X_t)_{t\in\RP}$,
given by $X_t := \sqrt{ S_t }$, $t \in \RP$, 
has the following asymptotic properties.
\begin{thmenumi}[label=(\alph*)]
\item\label{thm:limit-a} There is superdiffusive transience, namely, 
\begin{equation}
    \label{eq:escape-rate}
\lim_{t \to \infty} \frac{\log X_t}{\log t} = \frac{1+\gamma+\alpha}{2}, \as
\end{equation}
\item\label{thm:limit-b}
There is a distributional limit, namely, 
\begin{equation}
    \label{eq:limit-statement}
t^{-(1+\gamma+\alpha)/2}X_t \tod   \left(2\rho \int_0^1 s^\gamma \tY_s^\alpha \ud s\right)^{1/2},  \text{ as } t \to \infty,
\end{equation}
where $\tY$ in the limit follows~$\besq{\delta}{0}$.
\end{thmenumi}
\end{corollary}

\begin{remarks}
\phantomsection
\label{rems:limit}
\begin{myenumi}[label=(\roman*)]
 \item\label{rems:limit-i}
 Corollary~\ref{cor:bes-limit} is  direct from Theorem~\ref{thm:limit-new}.
Our proof of Theorem~\ref{thm:limit-new} establishes first a more general result, Theorem~\ref{thm:general} below,  identifying the  aspects of the exogenous noise process that are essential for the limiting behaviour in Theorem~\ref{thm:limit-new}. Differently put, Theorem~\ref{thm:general} extends the framework to admit a considerably larger class of processes: 
see \S\ref{sec:proof-overview} below for further remarks on the proofs.
We also note that the limit in~\eqref{eq:limit-statement2} implies an ``in probability'' version of~\eqref{eq:escape-rate2}, but not the a.s.~version. The a.s.~version  encapsulates the anomalous diffusion scaling
of the trajectory, whereas the ``in probability'' version 
is a marginal result that gives a looser reflection of this scaling. 
On the other hand, 
an ``in probability'' version of~\eqref{eq:escape-rate2} combined with a martingale decomposition and the self-similarity of SSCBI yield~\eqref{eq:limit-statement2}. 

 \item\label{rems:limit-ii} 
Under a mild additional
hypothesis, $X$ defined in Theorem~\ref{thm:limit-new} or 
Corollary~\ref{cor:bes-limit} 
satisfies dynamics~\eqref{eq:X_dynamics}, see Appendix~\ref{app:SDE_S}. 
It is more convenient 
to take as the primitive~$X^2=S$, because SDE~\eqref{eq:sq-X-dynamics}   fully specifies the process for a general continuous function $f$ (cf. Bessel process defined as a square-root of $\bessq$ in~\cite[Ch.~XI]{RevYor04}). 
\item\label{rems:limit-log}  
Since  $\gamma + \alpha \kappa >0$ by~\eqref{ass:parameters},  the scaling exponent
$(1+\gamma+\alpha\kappa)/2$ 
in~\eqref{eq:escape-rate2}
exceeds $1/2$, making the process~$X$ \emph{superdiffusive}. {In particular 
it follows from~\eqref{eq:escape-rate2} (here the a.s.~convergence is key) 
that $\lim_{t\to \infty} X_t =+\infty$, a.s., i.e., $X$ is transient.}  Moreover, by~\eqref{ass:parameters}, the values of $f$ on any compact set play no role.
\item\label{rems:limit-iv} 
The dimension parameter $\delta \in (0,\infty)$ of the squared-Bessel process $Y$  does not appear in the scaling exponent in~\eqref{eq:limit-statement}. The condition $\delta >0$
ensures that the Bessel process $\sqrt{Y}$ is \emph{diffusive},
which determines how the parameter~$\alpha$ in~\eqref{ass:parameters} enters the
scaling exponent in~\eqref{eq:escape-rate}--\eqref{eq:limit-statement}. In particular,  $X$ in Corollary~\ref{cor:bes-limit} exhibits superdiffusive transience for negative $\gamma\in(-\alpha,0)$,
even if $Y$ is topologically recurrent (i.e.~$\delta\in(0,2]$). A similar comment applies in the context of Theorem~\ref{thm:limit-new}.
\item\label{rems:limit-v} 
 We emphasize that the hypotheses in Theorem~\ref{thm:limit-new} 
 and Corollary~\ref{cor:bes-limit}
 permit
 arbitrary dependence between the process $Y$ and the driving Brownian motion $B$ in~\eqref{eq:sq-X-dynamics}. 
For example, in Corollary~\ref{cor:bes-limit}, since $Y$ with law
   $\besq{\delta}{y}$ is the unique strong solution of the stochastic differential equation 
 \begin{equation}
 \label{eq:sq_bessel}
 \ud Y_t = \delta \ud t + 2|Y_t|^{1/2}\ud W_t,  \text{ for } Y_0 = y\in\RP,
 \end{equation}
 for some Brownian motion $W=(W_t)_{t\in\RP}$, the limit in~\eqref{eq:limit-statement} requires only  that the Brownian motions $W$ and $B$ are adapted to the same filtration. In particular, $W$ and $B$   may be equal, independent or  have arbitrary stochastically evolving (adapted) covariation.
\item\label{rems:limit-vi} 
In the case $\alpha=0$, so that 
$f(t,y) \sim \rho t^\gamma$ as $t\to\infty$, the impact of 
 the process $Y$ 
on the large-scale dynamics of $X$ vanishes as $t\to\infty$.
In particular, the limit in Corollary~\ref{cor:bes-limit} is deterministic, and, in fact, the
convergence in distribution can in that case be strengthened to almost sure convergence. 
 Corollary~\ref{cor:bes-limit} can thus be viewed as a generalisation of certain results in~\cite{GO13} to non-polynomial time-inhomogeneous drift;
 see~\S\ref{sec:literature} below for some elaboration on this connection.
 \item\label{rems:limit-vii} 
In contrast to Remark~\ref{rems:limit}\ref{rems:limit-log},
 if we had $\gamma = \alpha =0$ (and hence $f$ asymptotically constant), 
then $X$ would be (essentially) a Bessel process of dimension $1 +2 \rho >1$.
In that case, the statement~\eqref{eq:limit-statement} is false, because 
the distributional limit of $t^{-1/2} X_t$ is \emph{random},
while the right-hand side of~\eqref{eq:limit-statement} is deterministic when $\alpha=0$.
This suggests that
the condition $\gamma + \alpha  >0$ in~\eqref{ass:parameters} cannot be omitted, essentially because Corollary~\ref{cor:bes-limit} requires the (exogenous) 
$Y$-driven stochastic drift 
to dominate the (endogenous) Brownian  noise in the evolution of~$X$.
\item \label{rems:limit-viii} We do not consider here the case where the process $Y$ in~\eqref{eq:X_dynamics}
 is \emph{ergodic}, because we anticipate rather different phenomena in that case. For example, if $f(s,y ) = h(y)$ depends only on~$y$, and if $h$ is integrable with respect to the stationary measure of~$Y$ on $\RP$, then it seems natural to suspect that $X$ behaves  similarly to a Bessel process with drift coefficient given via the stationary mean of~$h(Y)$. While interesting, this case seems unlikely to generate the anomalous diffusion that is the focus of the present work; see~\S\ref{sec:literature} below for more details.
\end{myenumi}
\end{remarks}

Remark~\ref{rems:limit}\ref{rems:limit-viii} suggests that for an ergodic process $Y$, the model $(S,Y)$
in~\eqref{eq:sq-X-dynamics} does not exhibit superdiffusivity. But what properties of $Y$ do guarantee anomalous diffusive behaviour in
 the model above? By Theorem~\ref{thm:general} below,   Assumption~\eqref{ass:noise} 
on the additive functional
\begin{equation}
    \label{eq:A-t-def}
   A_t := \int_1^t s^\gamma Y_s^\alpha \ud s , \text{ for } t \in [1,\infty), 
\end{equation}
ensures such behaviour.  Theorem~\ref{thm:limit-new} then follows by proving Assumption~\eqref{ass:noise} for {SSCBI} processes and identifying the weak limit.

\begin{description}
\item
[\namedlabel{ass:noise}{A$_Y$}]
Let $\alpha$ and $\gamma$ be as in~\eqref{ass:parameters}. Suppose $\kappa >0$ is such that the following hold. 
\begin{thmenumi}[label=(\alph*)]
\item\label{ass:noise-a} 
It holds that~$\int_{0}^{1} \Exp \bigl[ Y_t^\alpha \bigr] \ud t <\infty$.
\end{thmenumi}
Assume that $A_t$ in~\eqref{eq:A-t-def} satisfies the following hypotheses.
\begin{thmenumi}[label=(\alph*),resume]
\item\label{ass:noise-b} For a random variable $\tA$ in $\RP$,
$ A_t/t^{1+\gamma+\alpha \kappa} \tod \tA$ as $t\to\infty$.
\item\label{ass:noise-c} For every $\eps>0$, it holds that $\lim_{t\to\infty} t^\eps A_t/t^{1+\gamma+\alpha\kappa}  = \infty$, a.s.
\item\label{ass:noise-d}
It holds that $\sup_{t \in [1,\infty)} \Exp A_t/t^{1+\gamma+\alpha \kappa} <\infty$.
\end{thmenumi}
\end{description}

The next limit theorem requires no assumption beyond~\eqref{ass:noise} on the dynamics of $Y$.

\begin{theorem}
\label{thm:general}
Assume that~\eqref{ass:parameters}
holds with $\gamma + \alpha \kappa >0$, and that 
the adapted process $(S,Y,B)$
consists of 
an $\R$-valued Brownian motion $B$, 
a  process  $Y$  
on $\RP$ 
satisfying hypothesis~\eqref{ass:noise}, and an $\RP$-valued process 
$S$ satisfying SDE~\eqref{eq:sq-X-dynamics}.
Then, for $X_t = \sqrt{ S_t }$, it holds that
\begin{equation}
\label{eq:thm-general-a}
\lim_{t \to \infty} \frac{\log X_t}{\log t} = \frac{1+\gamma+\alpha\kappa}{2}, \as, 
\end{equation}
and, as $t \to \infty$, for $\tA$ the random variable in hypothesis~\eqref{ass:noise}\ref{ass:noise-b},
\begin{equation}
\label{eq:thm-general-b}
t^{-(1+\gamma+\alpha\kappa)/2}X_t \tod   \bigl( 2 \rho \tA \bigr)^{1/2}.
\end{equation}
\end{theorem}

 Theorem~\ref{thm:general} separates the analysis into two independent problems: a universal stochastic-kinetic argument and a process-specific verification of Assumption~\eqref{ass:noise}.
The proof of Theorem~\ref{thm:general} is given in~\S\ref{sec:limit-theory}. The proof of Theorem~\ref{thm:limit-new}
is presented in~\S\ref{sec:Y-lower-bounds},
where we establish that the process $Y$ with law $\sscbi{q}{c_0}{\delta}{y}$ satisfies Assumption~\eqref{ass:noise}.
 In \S\ref{sec:proof-overview} we give an overview of the proof strategy,
 after first (in~\S\ref{sec:motivation}) describing some motivation and relevant literature. 


We expect that other examples of exogenous noise processes satisfy hypothesis~\eqref{ass:noise} and hence can lead to further variations on Theorem~\ref{thm:limit-new}. Take, for example, $Y = |Z|$ where $Z$ is a symmetric $\nu$-stable process, $\nu \in (0,2)$. We do not pursue this example here.

\section{Motivation  and discussion}
\label{sec:literature}

\subsection{Motivation and literature}
\label{sec:motivation}

\subsubsection*{An example motivated by a self-interacting random walk.}
\begin{figure}[hbt]
    \centering
    \includegraphics[width=0.5\textwidth]{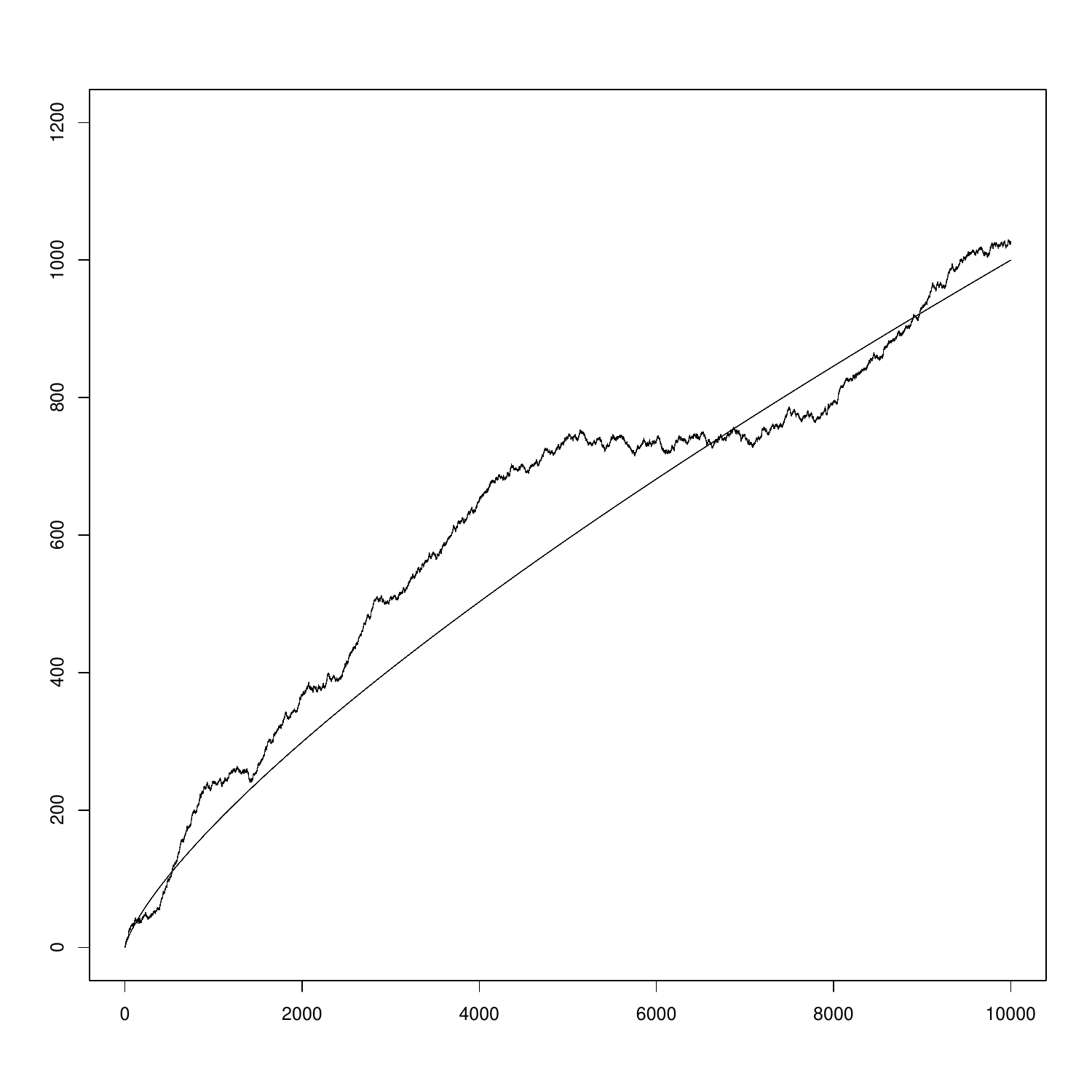} 
    \caption{Simulated realizations (Euler scheme, step size of~$1/10$) of a trajectory $t\mapsto X_t$ of $X$  
    (with $\gamma =0$, $\alpha =1/2$, i.e. $f(t,y)=\sqrt{y}$, and $\delta=1$)     and the graph $t \mapsto t^{p}$ for exponent $p = (1+\alpha+\gamma)/2 = 3/4$ on the time interval $[0,T]$ with $T = 10^5$.}
    \label{fig:3-4-trajectory}
\end{figure}

A discrete-time relative of the model in Corollary~\ref{cor:bes-limit} is studied in~\cite{cmsv},
motivated in part by a programme to study a certain self-interacting planar random walk.
The parameters in the present model that correspond, heuristically,
to the process described in~\cite{cmsv} are $\gamma =0$, $\alpha =1/2$, 
and $\delta=1$; informally for Brownian motions $B,W$ on $\R$ with arbitrary dependence,
$ \ud X_t = \rho | W_t |/X_t \ud t + \ud B_t$.
Note that $|W|$ (reflected BM) has the same law as $\sqrt{Y}$ where $Y$ follows $\besq{1}{0}$. 
By Corollary~\ref{cor:bes-limit}, as $t\to\infty$, $t^{-3/4} X_t$ converges weakly  to the law of
$(2\rho \int_0^1 | W_s| \ud s)^{1/2}$, which  can be expressed in terms of Airy functions (see~\cite[p.~349]{BSbook}). 
Figure~\ref{fig:3-4-trajectory} plots the graphs  of a path $t\mapsto X_t$ and $t\mapsto t^{3/4}$.

\subsubsection*{Stochastic kinetic dynamics.}
Various physical systems motivate \emph{stochastic kinetic models} for processes $(X,V) \in \R^d \times \R^d$
with dynamics of the form
\begin{equation}
    \label{eq:kfp}
 X_t = \int_0^t V_s \ud s , ~~~ V_t = \int_0^t G( V_s) \ud s + B_t ,\end{equation}
where $B$ is $d$-dimensional Brownian motion and $G : \R^d \to \R^d$. Here $V$ is an autonomous \emph{velocity process}
which feeds into the {\emph stochastic drift} of the $d$-dimensional process $X$.  Associated
to~\eqref{eq:kfp} is the so-called \emph{kinetic Fokker--Planck} equation (see~\cite{FT20,FT21}).
A classical example is Paul~Langevin's 1908 work on the \emph{confining}  case $G(v) = -v$ (see~\cite{langevin}),
 which was proposed to model physical Brownian motion.
In that case, $V$ is a $d$-dimensional Ornstein--Uhlenbeck process,
and, for some mean-zero Gaussian random variable $\bV$,
\begin{equation}
    \label{eq:langevin-velocity}
 V_t \tod \bV,  \text{ as } t \to \infty, \end{equation}
i.e., there is weak convergence to a stationary, isotropic velocity distribution.
Since $\bV$ is light-tailed and $V$ is rapidly mixing, then $X$ itself,
defined through the additive functional in~\eqref{eq:kfp}, 
will look like Brownian motion on large scales. 

Renewed interest in the system~\eqref{eq:kfp}, 
where $G(v)$ is of order $1/v$ for large $v$ (so that $V$ is a Bessel-like process), has been stimulated by both modelling of specific physical systems and realization that such processes can generate a variety of scaling behaviours corresponding to \emph{anomalous diffusion}; we refer to~\cite{bak,bg,kb,lp,FT20,FT21,BDL,GL}, and references therein. 

 {The mechanism for anomalous diffusion explored in these works is somewhat different from ours, as we explain. 
Indeed, these works maintain that $G$ be confining, taking $G (v) = -\delta/(1+v)$ (or  variations on this)
for $\delta >0$. For appropriate $\delta$, this leads to a \emph{heavy-tailed} stationary distribution $\bV$
for which~\eqref{eq:langevin-velocity} holds. 
In this case a competition between
the mixing rate of~$V$ and the index of the domain of attraction in which $\bV$ lives determines
the asymptotics of the additive functional~$X$ defined through~\eqref{eq:kfp}. 
 In contrast, our model $G$ is \emph{positive} and we are interested in
 velocity processes $V$ such that $V_t \to \infty$ in probability, as $t \to \infty$ (possibly recurrent, however). 
This gives a distinct mechanism for anomalous diffusion: rather than a stationary velocity process, we have a divergent in probability velocity process, but one whose self-similarity allows us to quantify the asymptotics of its integral.
Despite this distinct mechanism, the phenomena  are similar to those observed in e.g.~the $\delta \in (0,1)$ case of~\cite[Thm~2]{FT21}.}

A variation of~\eqref{eq:kfp} is to replace the dynamics of $X$ by
\begin{equation}
    \label{eq:kfp2}
 \ud X_t = \frac{f(t, Y_t)}{X_t} \ud t , \end{equation}
 for example, where $Y$ is an autonomous process. 
 Then
$\ud ( X_t^2) = 2 f(t,Y_t)  \ud t$, so, 
 in some sense, the generalization that~\eqref{eq:kfp2}
 provides to~\eqref{eq:kfp} is only in the form of the velocity process $V_t$.
 Our model~\eqref{eq:X_dynamics} extends~\eqref{eq:kfp2} further by introducing \emph{exogenous
 noise} into the dynamics of $X$.
 
\subsubsection*{Other stochastic drift models.}
From an entirely different direction, a number of two-dimensional models of a similar structure
to~\eqref{eq:kfp} and~\eqref{eq:kfp2} have been studied by Lefebvre and collaborators (see~\cite{lef2002a,lef2002b,lef2010,lef2012}),
motivated by lifetime studies and models of wear, among other applications. 
Specifically~\cite{lef2010,lef2012} formulate a version of~\eqref{eq:kfp2} for the model $(X,Y)$, where $\ud X_t = -   c (Y_t / X_t) \ud t$ and $Y$ is geometric Brownian motion. The focus of these papers is evaluating certain hitting distributions, rather than examining $t \to \infty$ asymptotics.

\subsubsection*{Long-time behaviour of stochastic interest-rate models.}
Models for the instantaneous interest rate (also known as the \textit{short rate}) $r$ typically use a mean-reverting diffusion as its stochastic representation of $r$. \cite{DEELSTRA1995163} (see also~\cite{Delbaen}) extend this framework to a model where the mean-reversion level of the diffusion $r$ is itself stochastic and has arbitrary dependence with the driving Brownian motion: the SDE for the short rate $r$ has constant negative mean-reversion rate $\beta<0$ and a stochastic mean-reversion level $Y$,
\begin{equation}
\label{eq:Y_ergodic}
\ud r_t = (\beta r_t+Y_t)\ud t+ \sqrt{r_t}\ud B_t,\qquad\text{with $\frac{1}{t}\int_0^tY_s \ud s\toas \obar{Y}$ as $t\to\infty$.}
\end{equation}
The main result in~\cite{DEELSTRA1995163} shows that $r$ is indeed mean-reverting, $\frac{1}{t}\int_0^tr_s\ud s\toas -\obar{Y}/\beta$ as $t\to\infty$, allowing, as in our Theorem~\ref{thm:general} above, arbitrary dependence between $Y$ and the Brownian motion $B$. In contrast to our setting, the process $\sqrt{r}$ in~\cite{DEELSTRA1995163} is diffusive.  Moreover, as long as $Y$ is ergodic as in~\eqref{eq:Y_ergodic} above, the process $r$ is likely to remain diffusive even if $\beta=0$, as it resembles the squared-Bessel process of positive dimension.

\subsubsection*{Time-inhomogeneous one-dimensional diffusions.}
In the case where the driving function~$f(t,y)$ in~\eqref{eq:X_dynamics} does not depend on~$y$ and is polynomial in $t$, the process $X$ falls into a class of one-dimensional diffusions with space-and-time-dependent drifts, studied in detail by~\cite{GO13} among a more general class of time-inhomogeneous diffusions on $\RP$, following earlier discrete-time work of~\cite{MV}; there are also structural links to the \emph{elephant random walk}~(\cite{bercu,bertoin}) and to the \emph{noise-reinforced Bessel process}~(\cite{bertoin-bessel}).
 
Specifically, a special case of dynamics~\eqref{eq:X_dynamics} is an SDE 
$\ud X_t = \rho t^\gamma X_t^{-1} \ud t + \ud B_t$, which is the  model of~\cite{GO13}
in the case $\alpha = -1$, $\beta = -\gamma$ (their parametrization). 
Here, for $\gamma >0$, \cite[Thm~4.10(i)]{GO13} states $\lim_{t \to \infty} t^{-\frac{1+\gamma}{2}} X_t = \left( 2\rho/(1+\gamma)\right)^{1/2}$, a.s., strengthening the log-scale asymptotic in~\eqref{eq:escape-rate} (where our $\alpha=0$) in this special case.

\subsection{Overview and discussion of the proofs}
\label{sec:proof-overview}

We outline our approach to Theorem~\ref{thm:general}. In what follows, the concrete case of Corollary~\ref{cor:bes-limit} can also be kept in mind,
where $Y$ is a squared-Bessel process.  
Recall that the adapted process $(S,Y,B)$, taking values in $\RP^2 \times \R$
satisfies the dynamics~\eqref{eq:sq-X-dynamics}.
In particular,
defining
\begin{equation}
  \label{eq:M-U-def}
 M_t :=  2 \int_0^{t} \sqrt{  S_s }  \ud B_s\quad \text{ and }\quad
  U_t :=  2 \int_0^{t} f(s,Y_s) \ud s,
\end{equation}
 we can write, for all $t \in \RP$,
\begin{equation}
  \label{eq:X2-decomposition}
S_t - S_0 =  U_t + t  +M_t . \end{equation}

The martingale decomposition~\eqref{eq:X2-decomposition} leads to a natural approach to the proofs,
in which the core step for both the weak limit~\eqref{eq:thm-general-b} and a.s.-asymptotic~\eqref{eq:thm-general-a} is to show that $M$ is a.s.~asymptotically negligible compared to~$U$. Hence $U$ determines the asymptotics of $X$, both in distribution  (via hypothesis~\ref{ass:noise}\ref{ass:noise-b}) and in 
a.s.-asymptotic sense (via hypotheses~\ref{ass:noise}\ref{ass:noise-c}--\ref{ass:noise-d}).

The growth rate of martingale $M$ is governed by its quadratic variation $\langle M\rangle= 4 \int_0^\cdot S_s \ud s$: by the Dambis--Dubins--Schwarz theorem,  $M = W_{\langle M\rangle}$ for a one-dimensional Brownian motion~$W$ time-changed by $\langle M\rangle$. A circular argument can be avoided, since~\eqref{eq:X2-decomposition} allows one to express $\Exp \langle M\rangle_t$ in terms of $\Exp U_t$, which is controlled by hypotheses~\ref{ass:noise}\ref{ass:noise-a} and~\ref{ass:noise-d}. An argument based on Doob's inequality then yields the desired a.s.~upper bound on $M$ (this is the content of Proposition~\ref{prop:M-martingale} below). 

While known results about Bessel processes give access to alternative reasoning if one is interested only in Corollary~\ref{cor:bes-limit}, 
we derive the Bessel case as a special case of the more general Theorem~\ref{thm:limit-new},
which shows that Theorem~\ref{thm:general} is applicable to a wide class of processes.
The hypotheses in assumption~\ref{ass:noise} are natural. Hypothesis~\ref{ass:noise}\ref{ass:noise-b}
is required for the distributional limit to exist.
We do not know any proof that works without some a.s.~upper bound, and hence without a
hypothesis like~\ref{ass:noise}\ref{ass:noise-d}. Combined with hypothesis~\ref{ass:noise}\ref{ass:noise-c} (which gives an a.s.~lower bound on $U$), this 
ensures that $U$ is genuinely dominant over~$M$ in~\eqref{eq:X2-decomposition}. It is the moments condition in~\ref{ass:noise}\ref{ass:noise-d} that imposes the restriction on $\alpha$ when $q < 1$ in Theorem~\ref{thm:limit-new}. 

We make some comments on the more apparently technical
conditions, which are to some extent necessary but where some variations are possible 
(say relaxing one at the expense of tightening another).
The finite-time integrability hypothesis~\ref{ass:noise}\ref{ass:noise-a} is mild (automatically satisfied in the Bessel case)
and provides initial control over~$\Exp \langle M\rangle_t$ that is not guaranteed by hypothesis~\ref{ass:noise}\ref{ass:noise-d}, since the integral in the definition 
of $A_t$ in~\eqref{eq:A-t-def} {starts at time~1}. In the integral in~\eqref{eq:A-t-def}, lower limit strictly greater than~$0$ avoids any possible singularity when $\gamma<0$; {so would using $1+s$ in place of $s$ in the integrand, at the expense of some  complications elsewhere (the expression being less amenable to scaling and self-similarity)}. Finally, we comment on the hypothesis on~$f$ given at~\eqref{eq:A}.
It is important in the present formulation that~\eqref{eq:A} provides an upper bound for $f(t,y)$ in terms of $y^\alpha$ for \emph{all} $t$ (see~\eqref{eq:f-upper-bound} below), since
it is used to reduce to the hypothesis~\ref{ass:noise}\ref{ass:noise-a}  what would otherwise be an hypothesis on~$\int_{0}^{1} \Exp \bigl[ f (t, Y_t) \bigr] \ud t$, which becomes a more indirect condition on~$f$. In the Bessel case, $Y$ spends almost all its time on order~$t$, and so it seems likely one could relax the hypothesis on~$f$ {to demand precise asymptotics only in a smaller region of $(t,y)$-space}, but we do not know any examples that would make  the extra work needed worthwhile.

\section{Strong law with stochastic drift}
\label{sec:limit-theory}

In this section we prove Theorem~\ref{thm:general} under Assumptions~\eqref{ass:parameters} and~\eqref{ass:noise}.
We start with a simple observation about $f$ under the hypothesis~\eqref{ass:parameters}.
By~\eqref{eq:A}, there exists $r_1>0$ such that  
 $f(t,y) \leq  (\rho+1) (1+t)^\gamma (1+ y)^\alpha$ 
 for all
 $(t,y) \in \RP^2$ with  $t + y\geq r_1$.
Since $f$ is continuous (and hence bounded on compacts),  
there exists $C_0 \in(0, \infty)$ such that
\begin{equation}
    \label{eq:f-upper-bound}
 f (t,y) \leq C_0(1 +  (1+t)^\gamma y^\alpha) \quad \text{ for all } (t,y) \in \RP^2. 
 \end{equation}
  By~\eqref{eq:f-upper-bound}, there exists $C < \infty$ such that $\sup_{0 \leq t \leq 1} f(t, y) \leq C (1+ y^\alpha )$. Hence the hypothesis that $\Exp \int_{0}^{1}  Y_t^\alpha \ud t < \infty$ from~\eqref{ass:noise} implies that $\Exp U_1  < \infty$,
  where $U$ is as defined in~\eqref{eq:M-U-def}. We use this fact, as well as the bound~\eqref{eq:f-upper-bound}, in several places in the proofs below.
 
Recall by~\eqref{ass:noise} that 
for $\alpha\in\RP$ and $ \gamma  \in (-\alpha, \infty)$, we have  $\beta\coloneqq \alpha+\gamma>0$.

\begin{proposition}
\label{prop:M-martingale}
     Assume~\eqref{ass:parameters} and \eqref{ass:noise}\ref{ass:noise-d} hold, and that $\alpha \in\RP$ and $\beta := \gamma + \alpha \kappa > 0$. 
     Then $M=(M_t)_{t \in \RP}$ in~\eqref{eq:M-U-def} is a martingale and there exists a constant $C \in(0, \infty)$ such that $\Exp ( M_t^2 ) \leq C t^{2+\beta}$ for all $t \in [1,\infty)$. Moreover, for any $\eps>0$, we have
     \begin{equation}\label{eq:M_t-upper-bound}
    \lim_{t \to \infty}  \sup_{0 \leq s \leq t} |M_s |/t^{1+(\beta/2) +\eps}  =0, \as
\end{equation}
\end{proposition}

\begin{proof}
 For any $N \in \RP$, set $\tau_N: = \inf\{t \in \RP \colon  S_t \geq N \}$. The quadratic variation of the local martingale 
 $M^{\tau_N}=(M_{t\wedge \tau_N})_{t\in\RP}$ is bounded: a.s., for all $t \in \RP$, 
 $\langle M\rangle_{t\wedge \tau_N}\leq 4\int_0^t S_{s\wedge \tau_N} \ud s\leq 4tN$.
 {Hence, by~\cite[Cor.~IV.1.25]{RevYor04}, 
 $M^{\tau_N}$ is a   martingale started at zero. } Thus $\Exp M_{t\wedge \tau_N}=0$ for all $t\in\RP$.
     Moreover,
 by definition~\eqref{eq:A-t-def}, the bound in~\eqref{eq:f-upper-bound} and~\eqref{ass:noise}\ref{ass:noise-d}, we have, for all $t \geq 1$,
 $0 \leq \Exp ( U_t - U_1 ) \leq 2C_0( t+\Exp A_t)\leq C_0' (1+t)^{1+\beta}$ for  a constant $C_0'\in \RP$.
Hence, since $U$ is non-decreasing and $\Exp U_1 < \infty$ (see the comment after~\eqref{eq:f-upper-bound} above), 
there exists  a constant $C_1\in \RP$ such that 
 \begin{align*}
\Exp S_{t\wedge \tau_N} 
& \leq \Exp S_0 + \Exp U_1 + \Exp ( U_t - U_1)  + t 
\leq C_1(1+t)^{1+\beta}, \text{ for all $t,N\in\RP$.}
 \end{align*}
 Hence $\Exp \langle M\rangle_{t\wedge \tau_N}\leq 4\int_0^t\Exp S_{s\wedge \tau_N}\ud s\leq 2C_1(1+t)^{2+\beta}$ for  $t,N\in\RP$. Since $\langle M\rangle_{t\wedge \tau_N}\uparrow \langle M\rangle_t$ as $N\to\infty$, 
 monotone convergence
 implies $\Exp\langle M\rangle_t\leq Ct^{2+\beta}$ for  $t\in[1,\infty)$  and  $C\coloneqq2^{3+\beta}C_1$. By~\cite[Cor.~IV.1.25]{RevYor04},
 $M$ is a martingale and $\Exp ( M_t^2 )=\Exp\langle M\rangle_t\leq Ct^{2+\beta}$ for $t\in[1,\infty)$.

The limit in~\eqref{eq:M_t-upper-bound} will follow by a Borel--Cantelli argument.
Doob's maximal inequality from~\cite[Thm~II.1.7]{RevYor04} and the $L^2$-bound 
 $\Exp M_t^2\leq Ct^{2+\beta}$ (for all $t\in[1,\infty)$)
yield
\begin{equation}\label{eq:Doob-ineq-M}
    \Pr \left( \sup_{0 \leq s \leq t} | M_s | >a \right) \leq \Exp M_t^2/ a^2 \leq C t^{2+\beta}/a^2 , \text{ for all  $t \geq 1$ and  $a >0$.}
\end{equation}
Fix arbitrary $\eps>0$. Set $t_n := 2^n$  and $a_n := t_n^{1+\beta/2 + \eps}$ for $n \in \ZP\coloneqq\N\cup\{0\}$. 
By the inequality in~\eqref{eq:Doob-ineq-M}, applied with $t=t_n$ and $a=a_n$, the probabilities of the events
\begin{equation*}
    E_n : = \left\{\sup_{0 \leq s \leq t_n} |M_s|>a_n \right\}\quad\text{are summable:}\quad
     \sum_{n = 0}^\infty   \bb{P}(E_n) \leq C \sum_{n = 0}^\infty t_n^{-2\eps}  < \infty.
\end{equation*}
Hence, by the 
Borel--Cantelli lemma,
there exists a (random)  $n_0 \in \bb{N}$, a.s., such that $E_n$ occurs for no $n \geq n_0$,
i.e.,
$\sup_{0 \leq s \leq t_n} | M_s | \leq a_n$ for all $n \geq n_0$.
Let $T := 2^{n_0}<\infty$ a.s. 
For every $t\in[1,\infty)$,
there exists a unique $k\in\ZP$, such that
$t_k = 2^k \leq t < 2^{k+1} =t_{k+1}$.
Thus,  for all $t \geq T$, we have
\begin{equation*}
\begin{aligned}
    \sup_{0 \leq s \leq t} |M_s| &\leq \sup_{0 \leq s \leq t_{k+1}} |M_s| \leq t_{k+1}^{1+\beta/2 +\eps} = 2^{1+\beta/2 + \eps} \big(2^{k}\big)^{1+\beta/2 + \eps}
\leq 2^{1+\beta/2 + \eps}  t^{1+\beta/2+\eps}.
\end{aligned}
\end{equation*}
Since $\eps>0$ was arbitrary, 
the limit in~\eqref{eq:M_t-upper-bound} holds almost surely.
\end{proof}

We now establish the following key result giving the almost-sure rate of escape for the additive functional $U$. In particular, by~\eqref{eq:X2-decomposition}, the almost-sure behaviour (as $t\to\infty$) of the process $S$ is dominated by  $U$.

\begin{proposition}
\label{prop:ratio-limit}
Suppose that $\alpha \geq 0$ and $\beta=\gamma + \alpha \kappa > 0$, and assume that~\eqref{ass:parameters}, \eqref{ass:noise}\ref{ass:noise-c}, and \eqref{ass:noise}\ref{ass:noise-d} hold.
Then, for every $\eps>0$,  the  process $U$ in~\eqref{eq:X2-decomposition} a.s.~satisfies
\begin{equation}\label{eq:U-log-limit}
t^{-\eps}< U_t /t^{1+\beta}< t^{\eps}  \text{ for all sufficiently large $t\in\RP$.}
\end{equation}
\end{proposition}
\begin{proof}
By~\eqref{eq:A} in~\eqref{ass:parameters}, there exists a large constant $r_\rho>1$ such that 
$f(t,y) 
\geq (\rho/2) t^\gamma y^\alpha $
for all $(t,y)\in[r_\rho,\infty)\times\RP$.
 By the definition of $U_t$ in~\eqref{eq:X2-decomposition},
    for all $t>r_\rho$,
    we have
    \begin{align*}
     \rho (A_t -  A_{r_{\rho}})=\rho \int_{r_{\rho}}^t s^\gamma Y_s^\alpha \ud s
        \leq    2\int_0^t f (s, Y_s ) \ud s =U_t, 
    \end{align*}
    where $A_t$ is defined in~\eqref{eq:A-t-def}.
Hypothesis~\eqref{ass:noise}\ref{ass:noise-c} on $A_t$ 
implies  $\lim_{t \to \infty} t^{\eps-(1+\beta)} A_t = \infty \as$ for every $\eps>0$, yielding the lower bound in~\eqref{eq:U-log-limit}. 

    To prove the upper bound in~\eqref{eq:U-log-limit}, note that the definition of $U_t$ in~\eqref{eq:X2-decomposition} and the upper bound on~$f$
    in~\eqref{eq:f-upper-bound} imply $U_t -U_1\leq 2C_0(t+A_t)$ a.s. Hence, by~\eqref{ass:noise}\ref{ass:noise-d}, for a constant
    $C_2 < \infty$,
    \[ 
     \Exp U_t = \Exp[U_t-U_1]+\Exp U_1\leq 2C_0(t+\Exp A_t)+\Exp U_1\leq C_2 t^{1+\beta}\quad\text{for all $t\in[1,\infty)$.}
    \]
    Pick any $\eps>0$. 
    The bound in the previous display  and Markov's inequality 
    yield
    \[ 
    \Pr (  U_t \geq t^{1+\beta+\eps} ) \leq \Exp  U_t/t^{1+\beta+\eps} \leq C_2 /t^{\eps}.\]   
Since $\sum_{n\in\ZP}1/t_n^{\eps}<\infty$, where $t_n\coloneqq 2^n$,
the Borel--Cantelli lemma implies 
$U_{t_n} \leq t_n^{1+\beta + \eps}$ 
for all but finitely many $n \in \ZP$, a.s.
As in the proof of Proposition~\ref{prop:M-martingale}, for every $t\in[1,\infty)$,
there exists a unique $k\in\ZP$, such that
$t_k = 2^k \leq t < 2^{k+1} =t_{k+1}$.
As $U$ is increasing (since $f \geq 0$)
\[   U_t  \leq  U_{t_{k+1}}    \leq (2^{k+1})^{1+\beta + \eps}= 2^{1+\beta+\eps} \cdot t_k^{1+\beta + \eps} \leq 2^{1+\beta+\eps} t^{1+\beta + \eps}.
\]
Since $\eps>0$ is arbitrary, the upper bound in~\eqref{eq:U-log-limit} holds for all large $t\in\RP$.
 \end{proof}
 
The next result establishes weak convergence of $U_t$, defined in~\eqref{eq:X2-decomposition} and suitably scaled,
to a random variable proportional to the limit $\tA$ in hypothesis~\eqref{ass:noise}\ref{ass:noise-b}.

\begin{lemma}
\label{lem:weak-limit-general}
Let $\alpha \geq 0$ and $\beta=\gamma + \alpha\kappa > 0$. Assume~\eqref{ass:parameters} and \eqref{ass:noise}\ref{ass:noise-b} hold.
 Then,  
 \begin{equation*}
 U_t/ t^{1+\beta} \tod 2 \rho \tA, \text{ as } t \to \infty.
\end{equation*}
\end{lemma}

\begin{proof}
Recall the definitions of $A_t$ in~\eqref{eq:A-t-def} and of 
$U_t$ in~\eqref{eq:X2-decomposition}.
Define also
\begin{align}
\label{eq:Dt}
D_{t}:= U_t/2 - \rho A_t =U_1/2+ \int_1^t f(s, Y_s) \ud s - \rho A_t\quad\text{ for any $t \in [1,\infty)$.}
\end{align}
Since $ U_t/t^{1+\beta}= 2(D_t+\rho A_t)/t^{1+\beta}$, by Slutsky's lemma
(see e.g.~\cite[p.~105]{Dur10}) and the weak limit in~\eqref{ass:noise}\ref{ass:noise-b}, it is sufficient to prove that $D_t/t^{1+\beta}$ converges to zero in probability. In particular, it suffices to show that for every $\eps', \eps''>0$ we have
\begin{equation}
\label{eq:D_t_conv_zero}
    \limsup_{t \to \infty} \Pr \bigl(  | D_t |/t^{1+\beta} > 2\eps' \bigr) \leq \eps''. 
\end{equation}

For every $\eps>0$, by~\eqref{eq:A}, there exists $r_\eps\in(1,\infty)$ such that, for every
 $s \in [r_\eps,\infty)$ and $y\in\RP$, we have  
 $\left| f(s,y)  - \rho s^\gamma y^\alpha \right|\leq 
 \eps s^\gamma (1+y)^\alpha$.
 Define
 $Z_\eps\coloneqq  U_1/2+\int_1^{r_\eps} f(s,Y_s) \ud s + \rho A_{r_\eps}$ and note that, by~\eqref{eq:Dt}, $| D_{r_\eps} | \leq Z_\eps$, a.s. Hence, for any $t>r_\eps$, by~\eqref{eq:A-t-def} and~\eqref{eq:Dt} again, we have 
\begin{align*}
| D_t | & \leq | D_{r_\eps} | + | D_t - D_{r_\eps} | \nonumber\\
& \leq Z_\eps + \int_{r_\eps}^t \left| f(s,Y_s)  - \rho s^\gamma Y_s^\alpha \right| \ud s  \leq Z_\eps + \eps \int_{r_\eps}^t s^\gamma (1+Y_s)^\alpha \ud s .
\end{align*}
Here $(1+Y_s)^\alpha = (1+Y_s)^\alpha \ind{Y_s \leq 1}
+ (1+Y_s)^\alpha\ind{Y_s >1} \leq 2^\alpha + 2^\alpha Y_s^\alpha$, so we get, for a constant $C < \infty$ depending on $\alpha$,
\begin{align}
\label{eq:E-Aut-1}
| D_t | \leq Z_\eps + C \eps t^{1+\gamma} + C \eps A_t.
\end{align}
Pick any $\eps',\eps''>0$. Fix small $\eps>0$, such that $C \eps t^{1+\gamma}/t^{1+\beta} < \eps'$ (which is possible since $\beta \geq \gamma$) and $\lim_{t \to \infty} \Pr ( \eps    A_t/t^{1+\beta} > \eps' ) = \Pr ( \tA > \eps'/\eps )$ (by~\eqref{ass:noise}\ref{ass:noise-b}, the limit holds for all but countably many $\eps$) and $\Pr ( \tA > \eps'/\eps ) < \eps''/2$.
Thus, by~\eqref{eq:E-Aut-1}, 
\begin{align*}
    \Pr \bigl(  | D_t |/t^{1+\beta} > 3\eps' \bigr) \leq \Pr \left(   Z_\eps/t^{1+\beta} > \eps' \right)
 + \Pr (  A_t/t^{1+\beta} > \eps'/\eps )\leq  \eps''/2+\eps''/2
\end{align*}
for all large $t\in\RP$.
Since 
$\eps', \eps''>0$ were arbitrary,
\eqref{eq:D_t_conv_zero} holds and
the lemma follows.
\end{proof}

We can now complete the proof of Theorem~\ref{thm:general}.

\begin{proof}[Proof of Theorem~\ref{thm:general}]
Let $\beta = \gamma +\alpha \kappa >0$.
Since 
 $S_t = S_0 + U_t + t+M_t$ and
$(S_0+t+M_t)/t^{1+\beta}\toas0$ as $t\to \infty$ by~\eqref{eq:M_t-upper-bound} of Prop.~\ref{prop:M-martingale}, the weak limit in Lemma~\ref{lem:weak-limit-general} and Slutsky's lemma (see~\cite[p.~105]{Dur10}) yield
$S_t/t^{1+\beta} \tod 2 \rho \tA$ as $t \to \infty$. Then~\eqref{eq:thm-general-b} follows by continuous mapping. 
By~\eqref{eq:U-log-limit} in  
Prop.~\ref{prop:ratio-limit},
for any $\eps>0$ a.s. 
$\log S_t/\log t= 1+\beta+\log\left((S_0+t+M_t)/t^{1+\beta}+U_t/t^{1+\beta}\right)/\log t\leq 1+\beta+\eps$
for all large $t$.
Similarly, for any $\eps\in(0,\beta/2)$, 
\eqref{eq:M_t-upper-bound} of Proposition~\ref{prop:M-martingale} and 
the lower bound in~\eqref{eq:U-log-limit} imply $1+\beta-\eps\leq \log S_t/\log t$ and  limit~\eqref{eq:thm-general-a} follows. 
\end{proof}

\section{Verifying the hypotheses of Theorem~\ref{thm:general}}
\label{sec:Y-lower-bounds}

The aim of this section is to verify Assumption~\eqref{ass:noise} for  the SSCBI $Y$.
This will, by Theorem~\ref{thm:general}, imply   Theorem~\ref{thm:limit-new}.
The focus in this section is the additive functional~$A_t$, defined in~\eqref{eq:A-t-def}. The next theorem collects the properties of $Y$ that we need.

\begin{theorem}
\label{thm:sscbi-integral-bounds}
Let $Y$ be a SSCBI with Laplace functional given at~\eqref{eq:sscbi} 
with parameters $\delta >0$ and $\kappa = 1/q \geq 1$, started at arbitrary $y\in\RP$.
Let $\alpha \in \RP$ and $\gamma \in \R$ with $1 +\gamma +\alpha\kappa >0$. Recall $A_t=\int_1^t s^\gamma Y_s^\alpha\ud s$, $t\in[1,\infty)$, defined in~\eqref{eq:A-t-def}. 
\begin{thmenumi}[label=(\alph*)]
\item\label{thm:sscbi-integral-bounds-a}
The following limit holds, 
\begin{equation}
    \label{eq:bessel-integral-log-limit}
\lim_{t \to \infty} \frac{ \log A_t}{\log t} =1+\gamma+\alpha\kappa, \as
  \end{equation} 
\item\label{thm:sscbi-integral-bounds-b}
Let $\tY$ follow~\eqref{eq:sscbi} with the same
parameters but started from $\tY_0 = 0$.
 Then, as $t \to \infty$,
 \begin{equation}
\label{lem:A-bessel-dist-limit}
  A_t / t^{1+\gamma+\alpha\kappa}\tod \int_0^1 s^\gamma \tY_s^\alpha \ud s,
\end{equation}
 where the right-hand side of~\eqref{lem:A-bessel-dist-limit} is an a.s.-finite random variable.
  \item\label{thm:sscbi-integral-bounds-c}
If (i) $q=1$, or (ii) $\alpha \in [0,q)$, then $\sup_{t \geq 1} \Exp  A_t/t^{1+\gamma+\alpha\kappa }  < \infty$.
  \end{thmenumi}
  \end{theorem}
  \begin{remark}
      \label{rem:sscbi-integral-bounds}
The condition $1+\gamma +\alpha \kappa >0$ in Theorem~\ref{thm:sscbi-integral-bounds} is weaker than the condition $\gamma+\alpha\kappa >0$ in Theorem~\ref{thm:limit-new}; the latter condition is needed to ensure the martingale part of the dynamics~\eqref{eq:sq-X-dynamics} is negligible compared to the kinetic drift.
  \end{remark}
 
Before giving the proof of Theorem~\ref{thm:sscbi-integral-bounds} we  introduce some notation. Recall that $Y = (Y_t)_{t\in\RP}$ 
with  $Y \sim \sscbi{q}{c_0}{\delta}{y}$  denotes a non-negative SSCBI process with Laplace transform given by~\eqref{eq:sscbi} started from $Y_0=y$.  Let $(L_t)_{t \in \RP}$ denote the Markov local time process of $Y$ at~$0$, and define the right-continuous inverse local time $\rho := (\rho_u)_{u \in \RP}$ by $\rho_u := \inf \{ t \in \RP : L_t > u\}$ for $u \geq 0$.
The next result contributes to the proof of Theorem~\ref{thm:sscbi-integral-bounds} by providing upper bounds on trajectories and expectation of $Y$.

\begin{lemma}
    \label{lem:sscbi-upper-bounds}
Suppose   $Y \sim \sscbi{q}{c_0}{\delta}{y}$ and $\kappa = 1/q$.
Then, for any $p \in [0,q)$, 
\begin{equation}
    \label{eq:sscbi-moment-bound}
    \Exp_y ( Y_t^p) < \infty \text{ for all } t \in \RP, \quad  \text{and} \quad
 \sup_{t \geq 1} t^{-p\kappa} \Exp_y (Y^p_t) < \infty.\end{equation}
 If $q=1$ then~\eqref{eq:sscbi-moment-bound} holds for all $p \geq 0$.
Moreover,  for all $q \in (0,1]$, with $\kappa = 1/q$, 
\begin{equation}
    \label{eq:sscbi-pathwise-bound}\limsup_{t \to \infty} \frac{\log Y_t}{\log t} \leq \kappa, \as \end{equation}
\end{lemma}
\begin{proof}
Using the identity (see 3.434 in~\cite[p.~333]{gr-tables}) valid for $p \in (0,1)$,
\[ x^p = \frac{p}{\Gamma (1-p)} \int_0^\infty \frac{1-\re^{-\lambda x}}{\lambda^{p+1}} \ud \lambda ,\]
we can write
\[ \Exp_y (Y_1^p) = \frac{p}{\Gamma (1-p)} \int_0^\infty \frac{1-\Exp_y ( \re^{-\lambda Y_1})}{\lambda^{p+1}} \ud \lambda .\]
Write $a :=  q c_0$ so that $1 - \Exp_y ( \re^{-\lambda Y_1}) = 1 - A(\lambda) B(\lambda)$ where $A(\lambda ) = (1+ a \lambda^q)^{-\delta} \in [0,1]$ and $B (\lambda ) = \exp \{ - \lambda y/ (1+ a\lambda^q)^{1/q} \} \in [0,1]$. Then
using the inequality $1-(1+z)^{-\delta} \leq \delta z$ for $z \in \RP$, we have
$1 - A(\lambda) \leq \delta a \lambda^q$, and using the inequality
$1-\re^{-z} \leq z$ for $z \in \RP$, we have
$1- B(\lambda) \leq \lambda y/ (1+ a\lambda^q)^{1/q} \leq \lambda y$. Therefore, for $0 < \lambda \leq 1$, 
\begin{align*} 1 - \Exp_y ( \re^{-\lambda Y_1}) & = 1 - A(\lambda) B(\lambda)
= 1 - A(\lambda) + A(\lambda) (1-B(\lambda)) \\
& \leq 1- A(\lambda) + 1 - B(\lambda) \\
& \leq \delta a \lambda^q + \lambda y \leq \lambda^q ( \delta q c_0 + y ).\end{align*}
Hence, for $p \in (0,q)$,
\begin{align*}
    \int_0^\infty \frac{1-\Exp_y ( \re^{-\lambda Y_1})}{\lambda^{p+1}} \ud \lambda
    & \leq(\delta q c_0 +y)  \int_0^1 \lambda^{q-p-1} \ud \lambda + \int_1^\infty \lambda^{-p-1} \ud \lambda \\
    & \leq \frac{\delta q c_0 +y}{q-p} +\frac{1}{p} < \infty.
\end{align*}
We conclude that, for $p \in (0,q)$,
\begin{equation}
    \label{eq:p-moment}
    \Exp_y (Y_1^p) \leq \frac{p}{\Gamma (1-p)} \left( \frac{\delta q c_0 +y}{q-p} +\frac{1}{p} \right).
\end{equation}
It follows from~\eqref{eq:p-moment} that, for every $p \in (0,q)$ and every $y_0 \in (0,\infty)$, $\sup_{y \in [0,y_0]} \Exp_y ( Y_1^p)<\infty$.
    Then the  scaling property says that $t^{-\kappa} Y_t$ under law $\Pr_y$ has the same law as $Y_1$ under $\Pr_{t^{-\kappa} y}$. 
    In particular, $\Exp_y ( Y_t^p) = t^{\kappa p} \Exp_{t^{-\kappa} y} (Y_1^p) < \infty$, for every $t >0$. Moreover,   
    $\sup_{t \geq 1} t^{-p \kappa}\Exp_y ( Y^p_t ) = \sup_{t \geq 1} \Exp_{t^{-\kappa} y} ( Y^p_1 ) <\infty$. This proves~\eqref{eq:sscbi-moment-bound} in the case $p \in(0,q)$, and the case $p=0$ is trivial. When $q=1$, $Y$ is a scalar multiple of a squared Bessel process, for which all moments exist and are locally bounded in the starting point. In particular, for any $p>0$, $\sup_{t\geq 1} t^{-p} \Exp_y ( Y_t^p) = \sup_{t \geq 1} \Exp_{y/t} (Y_1^p) \leq \sup_{0 \leq z \leq y} \Exp_z (Y_1^p) < \infty$, which implies~\eqref{eq:sscbi-moment-bound}.

To prove~\eqref{eq:sscbi-pathwise-bound} we show that for every $\eps>0$, $Y_t < t^{\kappa+\epsilon}$ for all $t$ sufficiently large, a.s. Write $v_t (\lambda) := \lambda (1 + q c_0 t \lambda^q)^{-1/q}$. Then from~\eqref{eq:sscbi} we have $\psi_t (\lambda ; y) := \Exp_y [ \re^{-\lambda Y_t}] = (1+qc_0 t \lambda^q)^{-\delta} \exp ( - y v_t (\lambda))$. Let $\tau_x := \inf \{ t \geq 0 : Y_t \geq x\}$. Let $T>0$. Then, on $\{ \tau_x \leq T \}$, by the strong Markov property at time $\tau_x$,
\[ \Exp_y [ 1 - \re^{-\lambda Y_T} \mid \cF_{\tau_x} ] = 1 - \psi_{T-\tau_x} (\lambda ; Y_{\tau_x} ) .\]
On $\{ \tau_x \leq T \}$, $Y_{\tau_x} \geq x$, so that  on $\{ \tau_x \leq T \}$, 
    \begin{align*} \psi_{T-\tau_x} (\lambda ; Y_{\tau_x} ) & \leq \psi_{T-\tau_x} (\lambda ; x ) \\
    & \leq \exp ( - x v_{T-\tau_x} (\lambda))  \leq \exp ( - x v_{T} (\lambda)) ,\end{align*}
    using the fact that $v_t (\lambda)$ is decreasing in $t$. Thus we obtain
\[ \Exp_y [ 1 - \re^{-\lambda Y_T} \mid \cF_{\tau_x} ] \geq 1 - \re^{-x v_T (\lambda)} , \text{ on } \{ \tau_x \leq T\} .\]
Consequently, a.s.,
\[ \ind { \tau_x \leq T } \leq \frac{\Exp_y [ 1 - \re^{-\lambda Y_T} \mid \cF_{\tau_x} ]}{ 1 - \re^{-x v_T (\lambda)} } .\]
Taking expectations we obtain
\begin{equation}
    \label{eq:maximal-inequality}
    \Pr_y \biggl( \sup_{0 \leq t \leq T} Y_t \geq x \biggr)
    = \Pr_y ( \tau_x \leq T ) \leq \frac{ 1 - \Exp_y [ \re^{-\lambda Y_T} ]}{ 1 - \re^{-x v_T (\lambda)} } .
\end{equation}
    For $0 < \lambda \leq 1$ we have the bound $ 1 - \Exp_y [ \re^{-\lambda Y_T} ] \leq \delta q c_0 T \lambda^q + y \lambda$. Thus from~\eqref{eq:maximal-inequality} we conclude, for all $\lambda \in (0,1]$, $T \in \RP$, and $x \in (0,\infty)$,
    \begin{equation}
    \label{eq:sscbi-maximal}
    \Pr_y \biggl(\sup_{0 \leq t \leq T} Y_t \geq x \biggr)
     \leq \frac{\delta q c_0 T \lambda^q + y\lambda}{ 1 - \re^{-x v_T (\lambda)} } .
\end{equation}
Now take $T_n = 2^n$, $x_n = T_n^{\kappa + \eps}$, and $\lambda_n = 1/x_n$.
Then
$x_n v_{T_n} (\lambda_n) = (1 + q c_0 T_n \lambda_n^q )^{-1/q}$, where
$T_n \lambda_n^q = T_n^{-q\eps}$, meaning that $x_n v_{T_n} (\lambda_n) \to 1$ as $n \to \infty$.
Hence, for all $n$ large enough, we get from~\eqref{eq:sscbi-maximal} that, for some constant $C<\infty$,
\[    \Pr_y \biggl( \sup_{0 \leq t \leq T_n} Y_t \geq x_n \biggr)
     \leq C \left(   T_n \lambda_n^q + y\lambda_n \right) ,\]
     which is summable in~$n$. Hence the Borel--Cantelli lemma implies that a.s., for all but finitely many $n\in \N$, $\sup_{0 \leq t \leq 2^n} Y_t < 2^{n (\kappa + \eps)}$.
     Every $T \geq 1$ has $2^n \leq T < 2^{n+1}$ for some $n \in \ZP$, so, a.s., for all $T$ large enough,
     \[ \sup_{0 \leq t \leq T} Y_t 
     \leq 
     \sup_{0 \leq t \leq 2^{n+1}} Y_t < 2^{(n+1) (\kappa + \eps)}
     \leq 2^{\kappa +\eps} \cdot T^{\kappa + \eps} ,
     \]
   which yields~\eqref{eq:sscbi-pathwise-bound}.
\end{proof}

The next result verifies the pathwise asymptotics needed for Theorem~\ref{thm:sscbi-integral-bounds}.

\begin{proposition}
\label{prop:sscbi-integral-bound}
Suppose   $Y \sim \sscbi{q}{c_0}{\delta}{y}$ and $\kappa = 1/q$.
Suppose that $\alpha\in \RP$ and $\gamma >-\alpha$. Then~\eqref{eq:bessel-integral-log-limit} holds.
\end{proposition}
\begin{proof}
First note that, by the pathwise upper bound from Lemma~\ref{lem:sscbi-upper-bounds}, for every $\eps>0$, 
$Y^\alpha_t < t^{\alpha \kappa+\eps}$ for all large $t$, so that, a.s.,
\[ A_t \leq \int_1^t s^{\gamma+ \alpha \kappa+\eps} \ud s = O ( t^{1+\gamma+\alpha\kappa+\eps} ) ,\]
and since $\eps>0$ was arbitrary, we get the upper bound half of~\eqref{eq:bessel-integral-log-limit}.
For the lower bound half of~\eqref{eq:bessel-integral-log-limit}, it suffices to prove that, for every $\eps >0$, a.s., for all $t$ sufficiently large,
\begin{equation}
    \label{eq:bessel-integral-bound-recurrent}
    A_t > t^{1+\gamma+\alpha \kappa - \eps} .
\end{equation}
    Suppose first that $ \delta \in (0,1)$ and $\alpha \in \RP$.
      We will first
    establish the case $\gamma = 0$ of~\eqref{eq:bessel-integral-bound-recurrent},
    then $\delta\geq 1$  and subsequently deduce the general case $\gamma \in (-\alpha,\infty)$.

    In the case $\delta \in (0,1)$, 
    the inverse local time (at $0$) process is a stable subordinator of index $\beta := 1-\delta$ (see~e.g.~\cite[p.~3]{mpu}), so $\Exp [ \re^{-\lambda \rho_u}] = \re^{-a u \lambda^\beta}$ for a constant $a \in (0,\infty)$.
    In particular,     $0$ is a recurrent point for $Y$, i.e.,
     $\inf \{ t \geq 0 : Y_t = 0\} < \infty$, a.s., so it suffices to suppose $Y_0=0$. 
     Then we have the self-similarity relation $( c^{-\kappa} Y_{ct} )_{t \in \RP} \eqd (Y_t)_{t \in \RP}$ for every $c>0$. 
     The process $Y$ and its inverse local time $\rho$ are jointly self-similar: for every $c>0$,
     \begin{equation}
         \label{eq:joint-ss}
         ( c^{-\kappa/\beta} Y_{c^{1/\beta} s} , c^{-1/\beta} \rho_{cu} )_{s,u \in \RP} \eqd (Y_s, \rho_u)_{s,u \in \RP}.
     \end{equation}
     
     Consider $A'_t := \int_0^t Y_s^\alpha \ud s$.
  We claim that the time-changed additive functional $A'_{\rho_u}$ is a subordinator. To see this, denote by $\cE$ the space of excursions of $Y$ away from $0$. Indexing excursions of $Y$ away from $0$ by local time $L$, It\^o excursion theory gives a Poisson point measure $\mathcal{N}$ on $\RP \times \cE$ with intensity $\ud \ell \nu ( \ud e)$ where $\ud \ell$ is Lebesgue measure and $\nu$ is the It\^o excursion measure. For $e \in \cE$ write $\zeta (e)$ for the duration of the excursion and define $I_\alpha (e) := \int_0^{\zeta(e)} e(s)^\alpha \ud s$. Then
  we have the representation
  \[ A'_{\rho_u} = \int_{(0,u] \times \cE} I_\alpha (e) \mathcal{N} ( \mathrm{d} \ell, \ud e ) .\]
  From this representation it follows that the non-decreasing process $(A'_{\rho_u})_{u \in \RP}$ has stationary and independent increments, so we verify that  $A'_{\rho_u}$ is a subordinator.

  Fix $c>0$. Then, using the change of variables $s = c^{1/\beta} v$ in the integral,
  \begin{align*}
      A'_{\rho_{cu}} & = \int_0^{\rho_{cu}} Y_s^\alpha \ud s  = c^{(1+\alpha \kappa)/\beta} \int_0^{c^{-1/\beta} \rho_{cu}} \bigl( c^{-\kappa/\beta} Y_{c^{1/\beta}v} \bigr)^\alpha \ud v ,
        \end{align*}
        and then the self-similarity property~\eqref{eq:joint-ss} shows that
        \begin{equation}
            \label{eq:additive-ss}
         ( A'_{\rho_{cu}} )_{u \in \RP}
        \eqd( c^{(1+\alpha \kappa)/\beta} A'_{\rho_u} )_{u \in \RP} .\end{equation}
        The self-similarity relation~\eqref{eq:additive-ss} identifies the subordinator $(A'_{\rho_{u}})_{u \in \RP}$ as a stable subordinator of index $\beta/(1+\alpha \kappa)$.

   By upper bound growth rates of subordinators (see~\cite[p.~92]{bertoin-book}) we have that, for every $\eps >0$, $\rho_t \leq t^{1/\beta+\eps}$ for all $t$ sufficiently large. On the other hand, by the corresponding lower bound (\cite[p.~92]{bertoin-book}, applicable since $\delta \in (0,1)$ and $\alpha \geq 0$), for every $\eps >0$, $A_{\rho_u} > u^{(1+\alpha \kappa)/\beta -\eps}$ for all $t$ sufficiently large. Combining these results and using monotonicity of $A_t$ shows that, for all $t$ sufficiently large,
    \[ A'_{t^{1/\beta+\eps}} \geq A'_{\rho_t} \geq t^{(1+\alpha\kappa)/\beta -\eps} ,\]
    which, after a change of variable, shows that, for any (deterministic or random) finite time $T$ and every $\eps>0$, 
   \begin{equation}\label{eq:gamma0-lb} 
\text{a.s.}\quad\int_T^t Y_s^\alpha \ud s \geq
     t^{1+\alpha\kappa -\eps}\quad\text{for all sufficiently large $t \in \RP$.}
   \end{equation}
In the case $\gamma =0$ and $\delta \in (0,1)$, the bound~\eqref{eq:gamma0-lb} 
directly implies~\eqref{eq:bessel-integral-bound-recurrent}.

If $\delta\geq1$, 
we have the representation $Y = \obar{Y}+ Y'$ where   $\obar{Y}$ and $Y'$ are independent with laws $\sscbi{q}{c_0}{\frac{1}{2}}{0}$ and $\sscbi{q}{c_0}{\delta -\frac{1}{2}}{y}$, respectively;
 this additivity in the immigration parameter  
  follows  from~\eqref{eq:sscbi}. Since $Y'_s\geq0$ for $s\in\RP$, by~\eqref{eq:gamma0-lb} (applied with $\delta=1/2$) we get 
\[ \int_0^t(\obar{Y}_s+Y'_s)^\alpha\ud s \geq \int_0^t\obar{Y}_s^\alpha\ud s\geq
     t^{1+\alpha\kappa -\eps} , \as, \] 
     for large $t$, implying that~\eqref{eq:gamma0-lb} holds for all $\delta>0$. This again implies the $\gamma =0$ case of~\eqref{eq:bessel-integral-bound-recurrent}. 
  
   Suppose $\gamma > 0$. By
   the pathwise upper bound in Lemma~\ref{lem:sscbi-upper-bounds}, 
   we have, for any $\eps \in (0,\gamma/\kappa)$, 
$\lim_{t \to \infty} Y_t^{1-\eps}/t^\kappa = 0$.
Hence, there exists (a random) $T\in\RP$, such that 
$t^\gamma \geq  Y_t^{\gamma/\kappa -\eps}$  and
$\int_0^t s^\gamma Y_s^\alpha \ud s \geq \int_T^t  Y_s^{\gamma/\kappa + \alpha -\eps} \ud s$
for $t\geq T$.
Thus, by~\eqref{eq:gamma0-lb} (with $\alpha+\gamma/\kappa -\eps$ in place of $\alpha$), \eqref{eq:bessel-integral-bound-recurrent} follows.
  If $\gamma \in(-\alpha,0)$,~\eqref{eq:bessel-integral-bound-recurrent} remains valid by~\eqref{eq:gamma0-lb} (with $T=1$) and 
   $\int_1^t s^\gamma Y_s^\alpha \ud s \geq t^\gamma \int_1^t   Y_s^\alpha \ud s$, concluding the proof.
   \end{proof}

Before being able to complete the proof of Theorem~\ref{thm:sscbi-integral-bounds}, we need one more technical result,  which provides a
 small-time bound for SSCBI started from $0$  needed for the convergence
in~\eqref{lem:A-bessel-dist-limit}.

\begin{lemma}
    \label{lem:small-times}
    Let $Y \sim \sscbi{q}{c_0}{\delta}{0}$. Then for every $\eps>0$, $Y_t = O (t^{\kappa-\eps})$ as $t \downarrow 0$.
\end{lemma}
\begin{proof}
    The starting point is the $y=0$ case of the maximal inequality~\eqref{eq:maximal-inequality}.
    Then by~\eqref{eq:sscbi}, $\Exp_0 [ \re^{-\lambda Y_T}] = (1+q c_0 T \lambda^q )^{-\delta}$.
    Using the inequality $1 - (1+z)^{-\delta} \leq C z$ for some constant $C<\infty$ and all $z \geq 0$, we get from~\eqref{eq:maximal-inequality} that
    \begin{equation}
        \label{eq:maximal-small}
        \Pr_0 \biggl( \sup_{0 \leq t \leq T} Y_t \geq x \biggr) \leq \frac{C T \lambda^q}{1 - \re^{-x v_T (\lambda)}} .
    \end{equation}
    Now fix $\eps \in (0,\kappa)$ and
    take $T_n = 2^{-n}$, $x_n = T_n^{\kappa-\eps}$, and $\lambda_n = 1/x_n$. Then $x_n v_{T_n} (\lambda_n) = (1+q c_0 T_n x_n^{-q} )^{-1/q} = (1+q c_0 T_n^{q\eps} )^{-1/q} \to 1$, as $n \to \infty$, so from~\eqref{eq:maximal-small} we get
    \[  \Pr_0 \biggl( \sup_{0 \leq t \leq T_n} Y_t \geq x_n \biggr) \leq C T_n^{q\eps} = C 2^{-nq\eps} . \]
    The Borel--Cantelli lemma then implies that, a.s., for all but finitely many $n \in \N$,
    $\sup_{0 \leq t \leq T_n} Y_t \leq 2^{-n(\kappa-\eps)}$. Every $T \in (0,1]$ has $2^{-n-1} < T \leq 2^{-n}$ for some $n \in \ZP$, so, for all $T>0$ sufficiently small,
    \[ \sup_{0 \leq t \leq T} Y_t \leq \sup_{0 \leq t \leq 2^{-n}} Y_t \leq 2^{-n(\kappa-\eps)}
    \leq 2^{\kappa-\eps} T^{\kappa-\eps} .\]
    The result follows.
\end{proof}



\begin{proof}[Proof of Theorem~\ref{thm:sscbi-integral-bounds}]
First, part~\ref{thm:sscbi-integral-bounds-a} of the theorem is proved in Proposition~\ref{prop:sscbi-integral-bound}, while part~\ref{thm:sscbi-integral-bounds-c} follows from~\eqref{eq:A-t-def} and the $ p= \alpha$ case of Lemma~\ref{lem:sscbi-upper-bounds}. It remains to prove the distributional convergence in part~\ref{thm:sscbi-integral-bounds-b}. The case where $\alpha =0$ is straightforward, so suppose $\alpha >0$.

By substituting $s=t v$ in definition~\eqref{eq:A-t-def} of $A_t$, we obtain 
\[
          A_t/  t^{1+\gamma+\alpha\kappa} = \int_{1/t}^{1} v^\gamma (t^{-\kappa}Y_{tv})^\alpha \ud v .\]
          We work with a convenient representation of $(t^{-\kappa}Y_{tv})_{v \in \RP}$.
          To this end, on the same probability space take independent $\tY = (\tY_v)_{v \in \RP}$ and  $Z = (Z_v)_{v \in \RP}$, where
          $\tY \sim \sscbi{q}{c_0}{\delta}{0}$ and $Z$ is the  continuous-state 
     branching process \emph{without immigration} that otherwise has the same parameters as~$Y$, i.e., $\Exp_y ( \re^{-\lambda Z_t} )$ is given by the right-hand side of~\eqref{eq:sscbi} 
     with $\delta =0$. Writing $v_t (\lambda ) = \lambda (1+q c_0 t \lambda^q)^{-1/q}$, this means that $\Exp_y ( \re^{-\lambda Z_t} ) = \re^{-y v_t (\lambda)}$. Write $X_v^{(t)} := \tY_v + t^{-\kappa} Z_{tv}$. Note that $( t^{-\kappa} Z_{tv} )_{v \in \RP}$ has the same law as $Z$ but started from $t^{-\kappa} y$. It follows that $(X^{(t)}_v)_{v \in \RP}$ has the same law as $(t^{-\kappa} Y_{tv} )_{v \in \RP}$. Hence to prove the distributional convergence in~\eqref{lem:A-bessel-dist-limit} it suffices to show 
     \begin{equation}
         \label{eq:coupled-convergence}
         \lim_{t \to \infty} \int_{1/t}^1 v^\gamma ( X^{(t)}_v )^\alpha \ud v = \int_0^1 v^\gamma \tY_v^\alpha \ud v, \as
     \end{equation}

Next consider $\tau := \inf \{ v \geq 0 : Z_v = 0\}$. Since $0$ is absorbing for $Z$, we have
\begin{align*}
    \Pr_y ( \tau > t ) = \Pr_y (Z_t >0) & = 1 -\lim_{\lambda \to \infty} \re^{- y v_t (\lambda)} \\
    & = 1 - \re^{-y (q c_0 t)^{-1/q}} .
\end{align*}
This last expression tends to $0$ as $t \to \infty$, so we conclude that $\Pr_y (\tau < \infty) = 1$. Set
\[ J_t := \int_{1/t}^1 v^\gamma ( t^{-\kappa} Z_{tv} )^\alpha \ud v.\]
Then, with the change of variable $s = tv$,
\[ J_t = t^{-1-\gamma-\alpha \kappa} \int_1^t s^\gamma Z_s^\alpha \ud s.\]
Here, since $\tau < \infty$ a.s., 
\[ \lim_{t \to \infty}  \int_1^t s^\gamma Z_s^\alpha \ud s
= \int_1^{1 \vee \tau} s^\gamma Z_s^\alpha \ud s < \infty , \as, \]
which means that $\lim_{t \to \infty} J_t = 0$, a.s., since $1+\gamma+\alpha \kappa >0$.

Now consider, for $\eta \in (0,1)$, 
\[ I_\eta = \int_\eta^1 v^\gamma \tY_v^\alpha \ud v ,\]
and note that $I_\eta < \infty$, a.s., for every $\eta \in (0,1)$, and $I_\eta \uparrow I_0$, a.s., as $\eta \downarrow 0$. We claim $I_0 < \infty$, a.s.; to see this it is enough to note that, by Lemma~\ref{lem:small-times}, as $\eta \downarrow 0$, for $\eps \in (0,1+\gamma+\alpha\kappa)$, 
\[ \int_0^\eta v^\gamma \tY_v^\alpha \ud v
\leq C \int_0^\eta v^{\gamma+\alpha \kappa -\eps} \ud v
= O ( \eta^{1+\gamma+\alpha \kappa -\eps} ) \to 0 .\]

Let $\eta \in (0,1)$. 
Then, since $X_v^{(t)} \geq \tY_v$, for all $t > 1/\eta$,
\[ \int_{1/t}^1 v^\gamma ( X_v^{(t)} )^\alpha \ud v \geq I_\eta ,\]
so taking $t \to \infty$ and then $\eta \downarrow 0$ we get
\begin{equation}
    \label{eq:liminf}
    \liminf_{t\to \infty} \int_{1/t}^1 v^\gamma ( X_v^{(t)} )^\alpha \ud v \geq I_0 , \as 
\end{equation}
On the other hand, for every $\alpha >0$ and $\theta >0$, there exists $C < \infty$ such that $(a+b)^\alpha \leq (1+\theta) a^\alpha + C b^\alpha$ for all $a, b \geq0$. Applying this inequality with $a = \tY_v$ and $b= t^{-\kappa} Z_{tv}$, we get
\begin{align*}
     \int_{1/t}^1 v^\gamma ( X_v^{(t)} )^\alpha \ud v 
     & \leq (1+\theta) I_{1/t} + C J_t \\
     & \leq  (1+\theta) I_{0} + C J_t.
\end{align*}
Recalling that $J_t \to 0$, a.s., as $t \to \infty$, we get
\begin{equation}
    \label{eq:limsup}
    \limsup_{t\to \infty} \int_{1/t}^1 v^\gamma ( X_v^{(t)} )^\alpha \ud v \leq (1+\theta) I_0 , \as 
\end{equation}
Since $\theta >0$ was arbitrary, we combine~\eqref{eq:liminf} and~\eqref{eq:limsup}
to verify~\eqref{eq:coupled-convergence}. This completes the proof of part~\ref{thm:sscbi-integral-bounds-b}.
\end{proof}



\begin{proof}[Proof of Theorem~\ref{thm:limit-new}]
Suppose that $Y$ follows $\sscbi{q}{c_0}{\delta}{y}$ , $\delta>0$, for a fixed $y \in \RP$.
We check the hypotheses of Theorem~\ref{thm:general},
i.e., we must check hypothesis~\eqref{ass:noise}.
First, Lemma~\ref{lem:sscbi-upper-bounds} shows that
hypothesis~\eqref{ass:noise}\ref{ass:noise-a} holds.
Theorem~\ref{thm:sscbi-integral-bounds}\ref{thm:sscbi-integral-bounds-a} implies hypothesis~\eqref{ass:noise}\ref{ass:noise-c},
Theorem~\ref{thm:sscbi-integral-bounds}\ref{thm:sscbi-integral-bounds-b} implies~\eqref{ass:noise}\ref{ass:noise-b}, with limit $\tA = \int_0^1 s^\gamma \tY^\alpha_s \ud s$, and Theorem~\ref{thm:sscbi-integral-bounds}\ref{thm:sscbi-integral-bounds-c} implies~\eqref{ass:noise}\ref{ass:noise-d}.
 Hence Theorem~\ref{thm:limit-new} is a consequence of Theorem~\ref{thm:general}.
\end{proof}

\appendix

\section{Existence, uniqueness and positivity for SDE~\texorpdfstring{\eqref{eq:sq-X-dynamics}}{}}
\label{app:SDE_S}
Let $(Y,B)$ be  an adapted process,
on a filtered probability space, 
where $B$ is a  scalar Brownian motion 
with respect to the given filtration
and $Y$ is a continuous adapted process in $\RP$. 

\subsubsection*{Existence and pathwise uniqueness.}
As noted in~\cite{Delbaen},
the SDE in~\eqref{eq:sq-X-dynamics} for $S$ is a special case of the Dol\'eans-Dade and Protter equation, 
driven by the 
semimartingale $B$ and an adapted continuous process $K=(K_t)_{t\in\RP}$, defined by $K_t\coloneqq \int_0^t(2f(s,Y_s)+1)\ud s$ for some continuous function $f:\RP\times\RP\to\RP$.
This equation is well-known~\cite{MR620991}
to possess a solution and
satisfy pathwise uniqueness (making every solution strong). In~\cite{Delbaen},  Deelstra \& Delbaen construct the solution $S$ of SDE~\eqref{eq:sq-X-dynamics} (which may \textit{a priori} take values in $\R$ with the volatility coefficient given by $\sqrt{|S_s|}$) from an Euler scheme approximation, and prove that if $S_0\in\RP$ then $S_t\in\RP$ for all $t\in\RP$ a.s. We stress that, beyond $Y$ and $B$ being adapted to the same filtration, no assumption is made on the dependence between the two processes (or on the dynamics of $Y$). 

\subsubsection*{Positivity.}
Consider now the strong solution $S$ of~\eqref{eq:sq-X-dynamics}, driven by $B$ and $K$ and started at $S_0\in\RP$. Assume in addition that  $1+2f\geq \delta_0$ a.s.\ for all $t>0$ and a constant $\delta_0>0$. Since the squared-Bessel SDE in~\eqref{eq:sq_bessel} has pathwise uniqueness, we may construct a squared-Bessel process $Z$, satisfying SDE~\eqref{eq:sq_bessel}  with ``dimension'' parameter $\delta_0$, started at $Z_0=0$ and driven by the Brownian motion $B$.
Note that by~\eqref{eq:sq-X-dynamics} and~\eqref{eq:sq_bessel} the quadratic variation of the semimartingale $S-Z$ equals
$\langle S-Z\rangle_t=4\int_0^t (\sqrt{S_s}-\sqrt{Z_s})^2\ud s$. In particular, since $(\sqrt{y}-\sqrt{x})^2\leq|x-y|$ for all $x,y\in\RP$, we have
$$
\int_0^t \ind{Z_s-S_s>0} \left|\sqrt{S_s}-\sqrt{Z_s}\right|^{-1}
\ud \langle S-Z\rangle_s\leq 4t<\infty \quad\text{a.s.\ for all $t\in\RP$.}
$$
By~\cite[Lem.~IX.3.3]{RevYor04}, the local time 
$L^0(S-Z)=0$ vanishes. Since $S_0\geq Z_0=0$, 
the Tanaka formula~\cite[Thm~VI.1.2]{RevYor04}
on  $[0,\tau_N]$, where $\tau_N\coloneqq \inf\{s\geq0:S_s\geq N\}$ 
for $N>0$, 
yields $\Exp \max\{0,Z_{t\wedge \tau_N}-S_{t\wedge \tau_N}\}\leq 0$ (since
 $s\delta_0\leq K_s$ for all $s\in\RP$ and $x\wedge y\coloneqq\min\{x,y\}$ for $x,y\in\R$).
Thus
$Z_{t\wedge \tau_N}\leq S_{t\wedge \tau_N}$, implying 
$S_t\geq Z_t$ a.s.\ for all $t\in\RP$ (since $\lim_{N\uparrow\infty}\tau_N=\infty$).

\subsubsection*{Bessel-type representation~\eqref{eq:X_dynamics} of $X=\sqrt{S}$.}
Since $f\geq0$, we can always take $\delta_0=1$, implying that the modulus of BM (i.e. $\sqrt{Z}$ for $\delta_0=1$) bounds $X=\sqrt{S}$ from below, making the point $0$ instantaneously reflecting for $X$. Moreover, if we had $f\geq \eps/2$ for  some $\eps>0$, a Bessel process $\sqrt{Z}$ with ``dimension''
$\delta_0=1+\eps$, which satisfies the Bessel SDE and, in particular, $\int_0^t\ud s/\sqrt{Z_s}<\infty$ for all $t\in\RP$, would bound $X$ from below. Thus,
\begin{equation}
\label{eq:finite_1_over_X_int}
\int_0^t (f(s,Y_s)/X_s) \ud s\leq \int_0^t (f(s,Y_s)/\sqrt{Z_s}) \ud s<\infty\quad\text{for all $t\in\RP$.}
\end{equation}
Once we know the integral in~\eqref{eq:finite_1_over_X_int} is finite, by considering the excursions of $X$ away from $0$, it is not hard to see that the quadratic variation at $t$ of the continuous process $X-\int_0^\cdot(f(s,Y_s)/X_s) \ud s$ equals $t$, making it, by L\'evy's characterisation, a Brownian motion and thus implying SDE~\eqref{eq:X_dynamics} for $X$.
However Assumption~\eqref{ass:parameters}
is not consistent with  $f\geq\eps/2$ when $\alpha>0$ (at large $t$, by~\eqref{eq:A}, $f(t,y)\to0$ as $y\downarrow0$). But we may  assume $f>0$ on $(0,\infty)\times(0,\infty)$, implying~\eqref{eq:finite_1_over_X_int} on the stochastic interval $[0,\tau^{(f)}_\eps)$,
where
$\tau^{(f)}_\eps\coloneqq\inf\{t\in\RP:f(t,Y_t)=\eps/2\}$ for small $\eps>0$.
Since, for 
$\delta\geq1$ and  $y>0$, the point $0$ is polar for 
$Y\sim\sscbi{q}{c_0}{\delta}{y}$, we have
$\tau^{(f)}_\eps\to\infty$ as $\eps\downarrow 0$ implying SDE~\eqref{eq:X_dynamics} for $X$ in this case. 

If $\delta\in(0,1)$, $Y$ hits zero a.s., bringing SDE~\eqref{eq:X_dynamics} into the realm of the   Bessel process with dimension~1 (which does not satisfy the corresponding SDE) at times $Y_t=0$. This case would require the analysis of the joint zeros $S_t=Y_t=0$, where the dependence in $(Y,B)$ clearly matters.

\section*{Acknowledgements}
\addcontentsline{toc}{section}{Acknowledgements}

The authors are grateful to Mikhail Menshikov and Vadim Shcherbakov for 
useful conversations. 
 CdC and AW are supported by EPSRC grant EP/W00657X/1. MB is supported by  EPSRC grant  EP/V009478/1. AM is  supported by
EPSRC grants  EP/V009478/1 and  EP/W006227/1 and, through  
The Alan Turing
Institute, by~EP/X03870X/1. Some of the work on this  paper was undertaken  during the programme ``Stochastic systems for anomalous diffusion'' (July--December 2024) hosted by the  Isaac Newton Institute, under EPSRC grant EP/Z000580/1.

\bibliographystyle{plain}
\bibliography{wp0}

\end{document}